\documentclass{article}

\usepackage{amssymb}
\usepackage{amsmath}

\def\thebibliography#1{\section*{References}\list
  {[\arabic{enumi}]}{\settowidth\labelwidth{[#1]}\leftmargin\labelwidth
    \advance\leftmargin\labelsep
    \usecounter{enumi}}
    \def\newblock{\hskip .11em plus .33em minus -.07em}
    \sloppy
    \sfcode`\.=1000\relax}
\newcommand{\refbook}[3]{{\sc #1}{\em\ #2}{\ #3}}
\newcommand{\refer}[5]{{\sc #1}{\ #2}{\em\ #3}{\bf\ #4}{\ #5}}

\newtheorem{lem}{Lemma}[section]

\newtheorem{teo}[lem]{Theorem}
\newtheorem{os}[lem]{Remark}
\newtheorem{defi}[lem]{Definition}
\newtheorem{prop}[lem]{Proposition}

\newcommand\sign{\mathop{\rm sign}}
\newcommand{\qed}{\thinspace\null\nobreak\hfill\hbox{\vbox{\kern-.2pt\hrule
 height.2pt depth.2pt\kern-.2pt\kern-.2pt \hbox to2.5mm{\kern-.2pt\vrule
 width.4pt \kern-.2pt\raise2.5mm\vbox to.2pt{}\lower0pt\vtop
 to.2pt{}\hfil\kern-.2pt \vrule
 width.4pt \kern-.2pt}\kern-.2pt\kern-.2pt\hrule height.2pt depth.2pt
 \kern-.2pt}}\par\medbreak}

\newcommand{\bound}{\mathcal{L}}
\newcommand{\A}{\mathcal{A}}
\newcommand{\R}{\mathbb{R}}
\newcommand{\K}{\mathbb{K}}
\newcommand{\C}{\mathbb{C}}

\newcommand{\Rp}{\textrm{\emph{Re}\,}}

\newcommand{\eps}{\varepsilon}

\newcommand{\Om}{\Omega}
\newcommand{\fra}{\mathfrak{a}}

\newcommand{\la}{\lambda}

\newcommand{\ov}{\overline}

\newcommand{\ds}{\displaystyle}

\date{}
\begin{document}

\title{Maximal regularity  for non-autonomous  Schr\"odinger  type  equations\thanks{Work supported by Universit\'{e} Bordeaux 1 and Universit\`{a} del Salento. The second author wishes to thank the Universit\`{a} del Salento for giving her the opportunity to visit the Universit\'{e} Bordeaux 1 and she is  grateful to the Institut de Math\'ematiques de Bordeaux, Universit\'{e} Bordeaux 1, for the hospitality.}}
\author{El Maati  Ouhabaz\thanks{Institut de Math\'ematiques de Bordeaux, Universit\'{e} Bordeaux 1,
351, cours de la Lib\'eration, F 33405 Talence cedex, France. e-mail: ElMaati.Ouhabaz@math.u-bordeaux1.fr }, Chiara  Spina\thanks{Dipartimento di Matematica ``Ennio De
Giorgi'', Universit\`a del Salento, C.P.193, 73100, Lecce, Italy.
e-mail: chiara.spina@unile.it, Chiara.Spina@math.u-bordeaux1.fr}}
\maketitle
\begin{abstract}
In this paper we  study  the maximal regularity property  for non-autonomous evolution equations  
$\partial_t u(t)+A(t)u(t)=f(t),\ u(0)=0.$
If the equation is considered on a Hilbert space $H$  and the operators  $A(t)$ are defined by sesquilinear forms $ \fra(t,.,.)$ we prove the maximal regularity   under a 
H\"older continuity assumption of $t \to \fra(t,.,.)$.   
In the non-Hilbert space situation we focus on Schr\"odinger type operators $A(t) := -\Delta + m(t, .)$ and prove $L^p-L^q$ estimates 
for a wide class of time and space  dependent  potentials  $m.$

\par

\end{abstract}

\section{Introduction}
Consider the Cauchy problem
\begin{equation} \label{pro-auto} 
\left\{
\begin{array}{ll}    
\partial_t u(t)+Au(t)=f(t),& t\in [0,T], \\
u(0)=0 \end{array}\right.
\end{equation}
where $\partial_t$ is the partial derivative with respect to the time variable.\\
The maximal regularity property for  autonomous evolution  equations  (the operator $A$ does not depend on $t$) of the type  (\ref{pro-auto})  has been widely investigated in the literature.  
One of the reasons why this  property is important lies in the fact that it allows to  study certain nonlinear problems. Indeed it is known that in this case the classical evolution operator approach fails whereas a technique based on the  maximal regularity for the linearized problem and the inverse function theorem allows to treat some quasilinear and fully nonlinear problems (see for example \cite{Pruss-Clement}). The maximal regularity for autonomous evolution equations  is  now well understood and   we refer to 
Section \ref{autonome} below for a more detailed description and references.  

For non-autonomous evolution equations (i.e., the operator $A$ depends also on the time variable $t$), the situation is much more difficult and it  is very less explored.  There are however 
several results in the litterature. 
Some authors have investigated  the case where the  operators $A(t)$ have the same domain (i.e. $D(A(t)) = D(A(0))$  is  independent  of $t \in [0, T]$). Pr\"uss and Schnaubelt \cite{Pruss-Schnaubelt}
 for example proved the maximal regularity under a time continuity assumption on $t \to  A(t)$. Arendt, Chill, Fornaro and Poupaud \cite{arendt} proved  a maximal  regularity result 
 by requiring that $D(A(t)) = D(A(0))$,  $t \to A(t)$ is  relatively continuous and  $A(t)$ satisfies the  maximal regularity for every fixed $t\in[0,T]$. Related results where also proved by 
 Amann \cite{amann}.\\
Concerning the case of operators whose domains depend on $t$,  Hieber and Monniaux \cite{hieber-monn} showed a maximal  regularity result in Hilbert spaces via the technique of pseudo-differential operators with operator-valued symbols provided that the family $\{A(t),\ t\in[0,T]\}$ satisfies the commutator Acquistapace-Terreni condition. Moreover they proved also \cite{hieber-monn2}  that a maximal $L^p-L^q$ regularity result holds by assuming  the Acquistapace-Terreni condition as well as some heat kernel bounds on $A(t).$\\ 
The Acquistapace-Terreni condition is rather strong, it requires a certain H\"{o}lder  regularity of $A$ (with respect to $t$) but it allows to study operators $A(t)$ with domains depending on $t$.\\
Such a condition consists in finding an estimate of the form
$$\|A(t)(\lambda-A(t))^{-1}(A(t)^{-1}-A(s)^{-1})\|_{\bound(H)}\leq c\frac{|t-s|^\beta}{1+|\lambda|^{1-\alpha}}$$ 
where $s,t\in[0,T]$ and $\lambda$ is in a sector contained in the resolvent set and which will be later specified. Evidently it is not easy to get a such estimate and it 
is one of the aims of this paper to  deal with  operators whose domains vary with $t$ without assuming the  Acquistapace-Terreni hypothesis. 
We will  prove a maximal regularity result for operators associated with sesquilinear 
 forms in a Hilbert space.  We will require that the forms have the same domain $V$ (but the domains of the operators may vary with $t$). 
Our idea lies in viewing each  operator $A(t)$ as an operator acting in the dual space $V'$ with domain 
 the whole space $V$. Using this we are in a somehow similar situation  than the case where the operators 
have the same domain. This allows one to obtain maximal regularity of the evolution equation but  considered on $V'.$ Under 
a H\"older continuity property of the sesquilinear forms (with respect to the time variable) we obtain the maximal regularity  of the evolution equation on $H.$ In order to achieve 
this we shall use in a crucial way   the previously  mentioned  results of  Hieber and Monniaux.
 We point out that the same restriction on the H\"{o}lderianity exponent appears also in some known results concerning the existence of classical solutions of (\ref{problem}) (see for example \cite[Section 5.4]{tanabe}).\\
Our maximal regularity result in  Hilbert spaces  applies to several examples, including uniformly  elliptic operators (with time dependent coefficients) on $L^2-$spaces  and 
also to  a  class of Schr\"{o}dinger operators with potentials  depending on the space and on the time variables. 
 In other words, we obtain $L^p-L^2$ a priori estimates for the corresponding parabolic equation.\\
The next step is to extend this and prove $L^p-L^q$ estimates for $q \not=2.$ For this, we shall concentrate on the case of 
Schr\"odinger type operators $\{-\Delta+m(t, \cdot),\ t\in[0,T]\}$ on $\R^N.$  Using $L^1-$estimates and domination arguments (Kato's  inequality)
 we  can prove  weak type $(1,1)$ 
boundedness of the operator $\partial_t(\partial_t-\Delta+m(t, \cdot))^{-1}$ and its adjoint.  From this, together with the $L^p-L^2$ estimates 
we obtain  that the family $\{-\Delta+m(t, \cdot),\ t\in[0,T]\}$ satisfies the $L^p-L^q(\R^N)$ maximal regularity   for all  $1<p,\ q<\infty$.\\
In the final section we explain how to prove  these $L^p-L^q$ estimates  for more general  operators  including elliptic operators on domains or
 Laplacians on some Riemannian manifolds.
 
 \bigskip
\noindent \textbf{Notation.} All over the paper we denote by $\rho(A)$ the resolvent set of $A$, by $\sigma(A)$ its spectrum, by $\bound(H)$ the space of all bounded operators on $H$ and by $\mu\{A\}$ the Lebesgue measure of a set $A$.

\section{Background material}
\subsection{The autonomous equations}\label{autonome}
Let $A$ be the generator of an analytic semigroup on some Banach space $X$ and consider  the autonomous parabolic problem associated with $A$
\begin{equation*}  
\left\{
\begin{array}{ll}    
\partial_t u(t)+Au(t)=f(t),& t\in [0,T], \\
u(0)=0. \end{array}\right.
\end{equation*}

As well known, in general, the derivative $\partial_t u$ of a solution of the above Cauchy problem is less regular than the right hand side $f$. We are concerned with the maximal regularity problem. Roughly speaking it means that we would like to know when this loss of regularity does not occur.

\begin{defi}
Given $p\in (1,\infty)$, we say that there is maximal $L^p$ regularity if for each $f\in L^p(0,T;X)$ there exists a unique $u\in W^{1,p}(0,T;X)\cap L^p(0,T; D(A))$ satisfying the above problem in the $L^p(0,T;X)$ sense.
\end{defi}

Many results concerning the maximal regularity property in the autonomous case can be found in literature. The Hilbert spaces case was first investigated. We mention some of them. For example De Simon \cite{Si64} proved a maximal $L^p(0,T;H)$ regularity result in Hilbert spaces $H$, then Sobolevskii \cite{Sob} proved that the maximal regularity property is independent of $p$. 
Da Prato and Grisvard \cite{DaGr} showed maximal  $L^p(0,T;X)$ regularity results in real interpolation spaces $X$. Lamberton \cite{Lam} proved that there is maximal $L^p$ regularity provided that $-A$ generates a bounded analytic semigroup on $L^2$ which acts as a contraction on $L^p$ for all $p\in [1,\infty]$. Dore and Venni \cite{DoVe} proved that there is maximal $L^p$ regularity provided that $A$ admits bounded imaginary powers (with an appropriate bound on their norms).  Weis \cite{weis} and \cite{weis1} gave necessary and sufficient conditions   for maximal regularity in terms of the so-called {\cal R}-boundedness of the resolvent or the semigroup of the operator $-A.$  Hieber and  Pr\"uss \cite{hieber-pruss}, Coulhon and Duong \cite{CoDu} obtained maximal $L^p$ regularity provided that the kernel of the semigroup generated by $A$ satisfies some gaussian upper bounds.

Let us recall why  the maximal regularity property implies some a-priori estimates for the solutions $u$ of the above evolution problem.\\ 
 In the  evolution equation $\partial_t$ is the distributional derivative with respect to $t.$ Note also that the
 operator $\partial_t: D(\partial_t)\rightarrow {L^p(0,T;X)}$, $1 \le p <\infty$ with domain  $D(\partial_t)=\{u\in W^{1,p}(0,T;X):\ u(0)=0\}$ is  a closed operator. It is even the generator of
 the translation semigroup (see for example \cite[Chapter II]{EngNag}). In the following we will deal also with the adjoint operator $\partial_t^*$ which acts  on 
 $D(\partial_t^*)=\{u\in W^{1,q}(0,T;X'):\ u(0)=u(T)=0\}$ where $q$ is such that $\ds \frac{1}{p}+\frac{1}{q}=1$ and so defined
$$\partial_t^*u=-\partial_t u$$
(see \cite[Section II.2]{gol} for a more detailed description of the adjoint operator).

\begin{os} 
If we consider $D(A)$ as a Banach space endowed with the graph norm, the maximal regularity property and the closed graph theorem imply that the operator $(\partial_t+A)^{-1}: L^p(0,T;X)\rightarrow L^p(0,T;D(\partial_t))\cap L^p(0,T; D(A))$ is continuous, hence there exists a positive constant $C$ such that for each $f\in L^p(0,T;X)$
$$\|(\partial_t+A)^{-1}f\|_{L^p(0,T;D(\partial_t))\cap L^p(0,T; D(A))}\leq C \|f\|_{L^p(0,T;X)}$$ or, if $u$ is a solution of the above problem,
$$\|u\|_{L^p(0,T;X)}+\|\partial_t u\|_{L^p(0,T;X)}+\|Au\|_{L^p(0,T;X)}\leq C\|(\partial_t+A)u\|_{L^p(0,T;X)}.$$
\end{os}

 \subsection{Non-autonomous equations}\label{n-autonome}
Suppose now that the operator $A$ is depending also on the time variable $t$ and consider  the non-autonomous parabolic problem associated with $A$ 
\begin{equation} \label{problem} 
\left\{
\begin{array}{ll}    
\partial_t u(t)+A(t)u(t)=f(t),& t\in [0,T], \\
u(0)=0. \end{array}\right.
\end{equation}

Generally, in this case  the domains $D(A(t))$ may vary with $t\in [0,T]$, therefore we cannot deduce maximal regularity results for (\ref{problem}) from the ones in the autonomous case by perturbation techniques.

\begin{defi}
Let $X$ be a Banach space and $p$ such that   $1<p<\infty$. The family $\{A(t),\ t\in[0,T]\}$ is said to belong to the class $MR(p,X)$ or equivalently we say that there is \emph{maximal $L^p-X$ regularity} for (\ref{problem}) if for each $f\in L^p(0,T;X)$ there exists a unique $u\in W^{1,p}(0,T;X)$  with $t\to A(t)u(t)\in L^p(0,T;X)$ satisfying (\ref{problem}) in the $L^p(0,T;X)$ sense.
\end{defi}

As in the autonomous case, the maximal regularity property and the closed graph theorem give some a-priori estimates for the solutions of (\ref{problem}) of the form
\begin{equation}  \label{apriori}
\|u\|_{L^p(0,T;X)}+\|\partial_t u\|_{L^p(0,T;X)}+\|Au\|_{L^p(0,T;X)}\leq C\|f\|_{L^p(0,T;X)}.
\end{equation}

As mentioned in the Introduction, a maximal regularity result in Hilbert spaces has been obtained  by M. Hieber and  S. Monniaux \cite[Theorem 3.2]{hieber-monn}. 
Since their result will be used in our proof, we  state it precisely.\\
The main assumptions  are the following  resolvent estimate H1) and the Acquistapace-Terreni condition H2). For $\theta\in (0,\pi)$, set $\Sigma_\theta:=\{z\in\C\setminus\{0\}; |\rm{arg}\ z|<\theta\}$.
\begin{itemize}
 \item [H1)] There exists $\theta\in \left(0,\ds\frac{\pi}{2}\right)$ such that $\sigma(A(t))\subset\Sigma_\theta$ for all $t\in [0,T]$ and for $\varphi\in (\theta,\pi)$ there exists $M>0$ such that 
$$\|(\lambda-A(t))^{-1}\|_{\bound(H)}\leq \frac{M}{1+|\lambda|},\ t\in[0,T],\ \lambda\in \C\setminus\Sigma_\varphi;$$
\item[H2)] there exist two constants $\alpha,\ \beta\in[0,1]$, $\alpha<\beta$, $\omega\in \left(\theta,\ds\frac{\pi}{2}\right),\ c>0$ such that
$$\|A(t)(\lambda-A(t))^{-1}(A(t)^{-1}-A(s)^{-1})\|_{\bound(H)}\leq c\frac{|t-s|^\beta}{1+|\lambda|^{1-\alpha}}$$
for $s,t\in[0,T],\ \lambda\in\C\setminus\Sigma_\omega$.
\end{itemize}

\begin{teo} [Hieber-Monniaux]      \label{Hie-Mon}
Let $1<p<\infty$, $T>0$ and $\{A(t),\ t\in[0,T]\}$ be a family of densely defined linear operators in a Hilbert space $H$ satisfying the  assumptions  H1) and H2).
Then $\{A(t),\ t\in[0,T]\}$ belongs to $MR(p,H)$.
\end{teo}

 \subsection{Forms and associated operators}\label{FO}
 
Let $H$ be a Hilbert space over $\K=\C$ or $\R$, $V$ be another Hilbert space contained in $H$ and  $\fra(t,\cdot,\cdot)$ a sesquilinear form defined on  $V\times V$ for every  fixed $t$ in a bounded interval $[0,T]$. We denote by $(\cdot,\cdot)$, $\Vert \cdot \Vert_H$ the inner product and the correspondent norm in $H$ respectively and by $((\cdot,\cdot))$, $\|\cdot\|_V$ the inner product and the correspondent norm in $V$. Assume  that
\begin{itemize}
\item[(i)] there exists a constant $M_0$ such that $\Vert u \Vert_H \leq M_0\|u\|_V$ for all $u\in V$;
\item[(ii)]  $\fra(t,\cdot,\cdot)$ is densely defined, i.e., $V$ is dense in $H$;
\item[(iii)] there exists a non-negative constant $M$ (independent  of $t$) such that 
$$  \vert \fra(t, u,v) \vert  \le M \|u\|_V \| v \|_V \, \,   \mbox{for all  }  u, v \in V \  \mbox{and  all }  t \in [0,T] ; $$
\item[(iv)] there exist a  positive number $\delta$ and a real number $k$ such that
$$\textrm{Re}\ \fra(t,u,u)\geq\delta\|u\|_V^2-k \Vert u \Vert_H^2 \ \mbox{ for all } t\in[0,T]  \mbox{ and all }    u\in V.$$
\end{itemize}
It is well known (see for example \cite[Chapter 1]{ouhabaz} or \cite[Chapter 2]{tanabe}) that, under these assumptions, one can associate with  $\fra(t,\cdot,\cdot)$ an  operator $A(t)$ defined by
\begin{eqnarray*}
& & D(A(t)) :=\{u\in V,  \, \exists\ g\in H:\ \fra(t,u,v)=(g,v)\ \forall\ v\in V\},\\ \nonumber
& & A(t)u =g.
\end{eqnarray*}
The operator $A(t)$  is a densely defined and accretive operator on $H.$ Its domain may depend on the variable $t.$ As mentioned in the Introduction, this latter fact is one of the main
difficulties  in order to obtain maximal regularity for the non-autonomous problem (\ref{problem}).\\ It possible to associate with $\fra(t,\cdot,\cdot)$ an operator $\A(t)$ whose domain is $V$ (hence independent of $t$) but it acts on a larger space than $H.$ 

Denote by $V'$ the dual (or anti-dual) space of $V$, that is the space of continuous linear (or ant-linear)  functionals $\phi$ on $V$.
Identifying  $H$ and its dual $H'$  yields
$$V\subset H\subset V'$$
 with continuous and dense imbedding.
We denote by $\langle\cdot,\cdot\rangle$ the dualization between $V$ and $V'$ (i.e. $\langle\phi,u\rangle$ denotes the value of $\phi$ at $u$ for $u\in V$ and $\phi\in V'$). 
In particular, if $\phi\in H$ and $u\in V$, then $\langle\phi,u\rangle=(\phi,u)$.
Fix $u\in V$ and consider the functional 
$$\phi(t,v):= \fra(t,u,v),\quad v\in V,\ t\in[0,T].$$
From the continuity assumptions (iii), it follows that the functional $\phi$ is continuous on $V$ and so it belongs to the dual space $V'$. By the Riesz representation theorem, there exists a unique vector $\A(t) 
u\in V'$ such that $\phi(t,v)=\langle\A(t) u,v\rangle$. The linearity of $\A(t)$  follows from the linearity of the form. By the continuity assumption  (iii) , we have
$$
 \|\A(t) u\|_{V'}=\sup_{\|v\|\leq 1}|\langle\A(t) u,v\rangle|=\sup_{\|v\|\leq 1}|\fra(t,u,v)|\leq M\|u\|_V.
$$
Thus $\A(t)$ is a continuous operator from $V$ into $V'$. Now let $A(t)$ be the operator associated with $\fra(t,\cdot,\cdot)$. By the density of $V$ in $H,$
 we see that $A(t)$ is the part of $\A(t)$ in $H$. This means that
$$D(A(t))=\{u\in D(\A(t)); \A(t) u\in H\}\quad\textrm{and}\ A(t)u=\A(t) u\quad \textrm{for}\ u\in D(A(t)).$$
For more information on $\A(t)$ (for each fixed $t$) see \cite[Section 1.4.2]{ouhabaz} or \cite[Section 2.2]{tanabe}. 

\section{Maximal regularity for non-autonomous equations in\\ Hilbert spaces} \label{sec3}
In this section,  $H$ and $V$ are Hilbert spaces such that $V \subset H$ with dense and continuous embedding. We use the same notation as in the previous section. \\We start with the following classical result due to J.L. Lions (see \cite{lions} or \cite{tanabe}). 
\begin{teo}\label{lions}
Let $\fra(t,\cdot,\cdot)$ be sesquilinear  forms  with the same domain  $V$ and satisfy  (i), (ii), (iii) and (iv) described in the previous section. Assume that 
$t \to \fra(t,u,v)$ is measurable for every $u, v \in V.$ For every $f \in L^2(0,T, V')$ there exists a unique $u \in W^{1, 2}(0,T, V') \cap L^2(0,T, V)$ satisfying
\begin{equation}\label{produal}
\partial_t u(t) + \A(t)u(t) = f(t), \ u(0) = 0.
\end{equation}
In other words,  the non-autonoumous Cauchy problem (\ref{produal}) has $L^2-V'$ maximal regularity. 
\end{teo}
If in addition $t \to \fra(t,u,v)$ is continuous, then (\ref{produal}) has maximal $L^p-V'$ regularity. 
\begin{prop}\label{pro} Let $\A(t), \, \fra(t,\cdot,\cdot)$ be as above and assume that 
for every $\varepsilon  > 0$, there exists a constant $\delta > 0$ such that
$$ |\fra(t,u,v) - \fra(s,u,v)| <  \varepsilon \|u\|_V\|v\|_V $$
for all $u, v \in V$ and all $t, s \in [0, T]$ such that $\vert t - s \vert < \delta.$ 
Then the problem (\ref{produal}) has $L^p-V'$ maximal regularity for every $p \in (1, \infty).$ 
\end{prop}
{\sc Proof.}   Fix $u\in V.$ We have 
\begin{eqnarray*}
\Vert \A(t) u - \A(s) u \Vert_{V'} &=& \sup_{\Vert v \Vert_V \le 1} \vert \fra(t, u, v) - \fra(s,u,v) \vert\\
& \le &   \sup_{\Vert v \Vert_V \le 1} \varepsilon \|u\|_V\|v\|_V \\
&=&  \varepsilon \|u\|_V.
\end{eqnarray*}
In particular, $t \to \A(t)$ is continuous from $[0, T]$ to ${\cal L}(V, V')$. The conclusion of the proposition follows then from \cite{Pruss-Schnaubelt}. \qed

Note that by using \cite[Theorem 2.7]{arendt} we  can weaken  the continuity assumption  in the previous proposition.

We are in position to state  the maximal regularity property for  operators $A(t)$ on $H.$ 
\begin{teo}  \label{monn}
 Let $\fra(t,\cdot,\cdot)$ be sesquilinear  forms  with the same domain  $V$ and satisfy  (i), (ii), (iii) and (iv) described above. Suppose moreover that
$\fra(t,u,v)$ is H\"{o}lder continuous in $t$ in the following sense: there exist $\beta > \frac{1}{2}$, $K>0$ such that
\begin{equation} \label{Holder}
|\fra(t,u,v) - \fra(s,u,v)|\leq K|t-s|^\beta\|u\|_V\|v\|_V
\end{equation}
for all $s,\ t\in[0,T]$ and $u,\ v\in V$. Let $A(t)$ be the operator associated with the form $\fra(t,\cdot,\cdot)$ as previously defined. Then $\{A(t),\ t\in[0,T]\}$ belongs to $MR(p,H)$. In other words, (\ref{problem}) has $L^p-H$ maximal regularity. 
\end{teo}

Before we prove this result, let us mention that $ \{A(t),\ t\in[0,T]\}$ has $L^p-H$ maximal regularity {\it iff} for every  $\delta > 0$, $ \{A(t) + \delta I,\ t\in[0,T]\}$ has $L^p-H$ maximal regularity. Therefore, we may assume without loss of generality that assumption (iv) holds with $k = 0.$ 

We first prove  a preliminary lemma.
\begin{lem}  \label{risolv}
Let $A(t)$ be the operator associated with $\fra(t,\cdot,\cdot)$ under the assumptions above (with $k = 0$ in (iv)). Then there exists a positive constant $C$ such that for all $f\in H$, $t\in[0,T],$ $\lambda\in\rho(A(t))$, $\Rp\lambda\leq 0$,   
$$\|(\lambda-A(t))^{-1}f\|_V \leq C|\lambda|^{-\frac{1}{2}} \Vert f \Vert_H.$$
\end{lem}
{\sc Proof.} Let $f\in H$, $\lambda\in\rho(A(t))$, $\Rp\lambda\leq 0$ and set  $u=(\lambda-A(t))^{-1}f$.  Then, by the definition of $A(t)$,
$$(f,v)=\lambda(u,v)-\fra(t,u,v)$$
for every $v\in V$. If we take $v=u$, the previous becomes
\begin{equation}  \label{eq0}
(f,u)=\lambda(u,u)-\fra(t,u,u)
\end{equation}
from which it follows that 
\begin{equation} \label{eq1}
\|u\|_H\|f\|_H\geq \Rp [-(f,u)]= \Rp \fra(t,u,u)-\Rp \lambda \|u\|_H^2\geq \delta\|u\|8V^2
\end{equation}
and hence, by (\ref{eq0}) again and the continuity assumption on the form, it follows that
$$|\lambda|\|u\|_H^2\leq \|u\|_H\|f\|_H+M\|u\|_V^2\leq \left(1+\frac{M}{\delta}\right)\|f\|_H\|u\|_H$$
where $M$ is the constant in the continuity assumption of the form. Therefore we deduce
\begin{equation} \label{eq2}
|\lambda|\|u\|_H\leq\left(1+\frac{M}{\delta}\right)\|f\|_H.
\end{equation}
We conclude by observing that, by (\ref{eq1}) and (\ref{eq2})
$$\delta\|u\|_V^2\leq \|f\|_H\|u\|_H\leq \frac{1}{|\lambda|}\left(1+\frac{M}{\delta}\right)\|f\|_H^2$$
which is the claim.\qed  

{\sc Proof of Theorem \ref{monn}}. Since the family $\{A(t),\ t\in[0,T]\}$ belongs to $MR(p,H)$ if and only if this is true for $\{A(t)+\delta I,\ t\in[0,T]\}$ with $\nu$ arbitrary constant, we can suppose that $A(t)$ is invertible for all $t\in [0,T]$.
We would like to apply the Theorem \ref{Hie-Mon}, so we need to verify that the assumptions H1) and H2) are satisfied. 
Let us first observe that the operator $-A(t)$ generates an analytic semigroup both in $H$ and in $V'$ (see \cite[Chapter 1]{ouhabaz}),  therefore the assumption H1) is verified for some $\theta\in(0,\frac{\pi}{2})$. Let us turn our attention to the second assumption. Denote by $L$ the operator $A(t)(\lambda-A(t))^{-1}(A(t)^{-1}-A(s)^{-1})$ defined on $H$ and by $\mathcal{L}$ the analogous operator with $\A$ instead of $A$ defined on $V'$. 
Let $\omega\in \left(\theta,\ds\frac{\pi}{2}\right)$,   
$\lambda\in\Sigma_\omega$, $u,\ v\in V$, $t,s\in [0,T]$. Since the domain of $\A(t)$ coincides with $V$ and so it is independent of $T$, then $\A(t)$ and $(\lambda-\A(t))^{-1}$ can commute in the expression of $\mathcal{L}$ and we can write
\begin{eqnarray*}
 \langle\mathcal L u,v\rangle &=&\langle \A(t)(\lambda-\A(t))^{-1}(\A(t)^{-1}-\A(s)^{-1})u,v\rangle\\&= &
\langle (\lambda-\A(t))^{-1}\A(t)(\A(t)^{-1}-\A(s)^{-1})u,v\rangle\\&=&
\langle (\lambda-\A(t))^{-1}(\A(s)-\A(t))\A(s)^{-1}u,v\rangle\\&=&
\langle(\A(s)-\A(t))\A(s)^{-1}u,(\ov{\lambda}-\A^*(t))^{-1}v\rangle\\&=&
\langle\A(s)\A(s)^{-1}u,(\ov{\lambda}-\A^*(t))^{-1}v\rangle-\langle\A(t)\A(s)^{-1}u,(\ov{\lambda}-\A^*(t))^{-1}v\rangle
\end{eqnarray*}
where $\A^*$ is the adjoint operator of $\A$.
By the definition of $\A$,
$$\langle\mathcal L u,v\rangle=a(s,\A(s)^{-1}u,(\ov{\lambda}-\A^*(t))^{-1}v)- \fra(t,\A(s)^{-1}u,(\ov{\lambda}-\A^*(t))^{-1}v).$$
The H\"{o}lderianity assumption on the form implies that
$$|\langle\mathcal L u,v\rangle|\leq K|t-s|^\beta\|\A(s)^{-1}u\|_V\|(\ov{\lambda}-\A^*(t))^{-1}v\|_V$$
for all $u,\ v\in V$.
By  Lemma \ref{risolv},
$$\|(\ov{\lambda}-\A^*(t))^{-1}v\|_V\leq C|\lambda|^{-\frac{1}{2}}\|v\|_H $$
with $C$ independent of $t$ and, analogously,
$$\|\A(s)^{-1}u\|_V\leq C\|u\|_H$$ for some other positive constant $C$.
We obtained that 
$$|\langle \mathcal{L} u,v\rangle|\leq C\frac{|t-s|^\beta}{|\lambda|^\frac{1}{2}}\|u\|_H\|v\|_H.$$
Hence $\mathcal L u\in H$ for $u\in V$ and, by the density of $V$  in $H$, 
$$\|L\|_{\bound({H})}\leq C\frac{|t-s|^\beta}{|\lambda|^\frac{1}{2}}.$$
This proves that assumption H2) is satisfied with $\alpha=\frac{1}{2}$ and $\beta$ as in the H\"{o}lderianity assumption of the form. By Theorem \ref{Hie-Mon} we deduce the claim.\qed

\vspace{.5cm}
\noindent{\bf Examples.}\\

1)  {\it Uniformly elliptic operators on domains.} Let $\Omega$ be an open subset of $\R^N$ ($N \ge 1$) endowed with the Lebesgue measure $dx$. Denote by $V$ a fixed closed subset of the Sobolev space $H^1(\Omega)$ which contains the space of $C^\infty-$functions with compact supports in $\Omega.$  We define on 
$H := L^2(\Om, dx)$ the sesquilnear forms (here $t \in [0,T]$ and $T > 0$ is fixed)
\begin{eqnarray*}
\fra(t,u,v) &=& \sum_{k,j=1}^N \int_\Omega a_{kj}(t,x) \partial_j u\cdot\partial_k v\ dx + \sum_{k=1}^N \int_\Omega \left[ a_k(t,x) \partial_k u\cdot v + b_k(t,x) u\cdot \partial_k v \right] dx\\
&&  + \int_\Omega m(t,x) u\cdot v\ dx\\
D(\fra(t,\cdot,\cdot)) &=& V \,  \ \mbox{for all } \ t \in [0,T].
\end{eqnarray*}
Here we assume that the coefficients $a_{kj}, a_k, b_k$ and $m$ satisfy the following  conditions  (in which $\eta$ and $M$ are positive constants independent of $t$). 
\begin{equation}\label{a1}
\max \{ \vert a_{kj}(t,\cdot) \vert, \vert a_k(t,\cdot) \vert, \vert b_k(t,\cdot) \vert, \vert m(t,\cdot) \vert \} \le M \ (a.e.\ x \in \Omega),
 \end{equation}
\begin{equation}\label{a2}
 \sum_{k,j=1}^N a_{kj}(t,x) \xi_k \xi_j \ge \eta \vert \xi \vert^2 \ \mbox{for all} \, \xi =(\xi_1,...,\xi_N) \in \R^N, 
 \end{equation}
\begin{equation}\label{a3}
 \max\{  \vert a_{kj}(t,x) - a_{kj} (s,x) \vert, \vert a_{k}(t,x) - a_{k} (s,x) \vert, \vert b_{k}(t,x) - b_{k} (s,x) \vert, \vert m(t,x) - m(s,x) \vert \} \le M \vert t - s \vert^\beta.
 \end{equation}
Here $\partial_k := \frac{\partial}{\partial x_k}$ and the estimates  hold for all $t, s \in [0,T]$ and $a.e. x \in \Omega.$ \\
The first  two  assumptions (\ref{a1}) and (\ref{a2}) imply that the forms $\fra(t,\cdot,\cdot)$ are closed and  satisfy the assumptions (i), (ii), (iii) and (iv) of Theorem \ref{monn}. 
Therefore, if   (\ref{a3})  holds with some constant $\beta > 1/2$ then the problem (\ref{problem}) has $L^p-L^2(\Omega)$ maximal regularity. The operators $A(t)$ are now the associated 
operators with the forms $\fra(t,\cdot,\cdot)$. These are time-dependent uniformly elliptic operators subject to the boundary conditions fixed by $V.$ \\

 2) {\it Schr\"odinger operators.}   We consider a particular form of the operators of the previous example. However, we want to include the case where the potential 
 $m$ is not bounded. We concentrate on $A(t) := -\Delta + m(t,\cdot)$ on $L^2(\R^N)$ but we could consider similar operators on domains with general boundary conditions. \\
 We define 
 the forms
 \begin{eqnarray*}
\fra(t,u,v) &=& \sum_{k=1}^N \int_{\R^N}  \partial_k u \cdot \partial_k v\ dx  + \int_{\R^N} m(t,x) u\cdot v\ dx\\
D(\fra(t,\cdot,\cdot)) &=& \{ u \in H^1(\R^N), \int_{\R^N} m(t,x) \vert u \vert^2 dx < \infty \}. 
\end{eqnarray*}
We assume that there exists a non-negative potential $W \in L^1_{loc}(\R^N, dx)$ such that $m$ satisfies the following properties 
(in which $c_1, c_2$    are positive constants and $\beta > 1/2$). 
\begin{equation}\label{a4}
c_1 W(x) \le m(t,x) \le c_2 W(x) \ (a.e.\  x \in \R^N) \ \mbox{and all } \ t \in [0,T], 
\end{equation}
\begin{equation}\label{a5}
 \vert m(t,x) - m(s,x) \vert  \le c_2 W(x) \vert t - s \vert^\beta  \ (a.e.\ x \in \R^N) \ \mbox{and all } \ t, s \in [0,T].
\end{equation}
Under these assumptions, it is clear that 
$$D(\fra(t,\cdot,\cdot)) = V := \left\{ u \in H^1(\R^N), \ \int_{\R^N} W(x) \vert u(x) \vert^2 dx < \infty \right\} \ \mbox{for all } \ t \in [0,T].$$
The space $V$ endowed with  the norm
$$\Vert u \Vert_V := \left[  \int_{\R^N} \vert \nabla u \vert^2 dx +  \int_{\R^N} \vert  u \vert^2 dx +  \int_{\R^N} W \vert  u \vert^2 dx \right]^{1/2}$$
is a Hilbert space. The form $\fra(t,\cdot,\cdot)$ is closed and its  associated operator  is the Schrodinger operator
$ A(t) := - \Delta + m(t,\cdot).$\\
Under (\ref{a4}) and (\ref{a5}) the forms  satisfy the assumption of Theorem \ref{monn}  if  $\beta >  1/2.$ Therefore, the time dependent Schr\"odinger equation
\begin{equation} \label{proSch} 
\left\{
\begin{array}{ll}    
\partial_t u(t) - \Delta u(t) + m(t,\cdot) u(t) =f(t),& t\in [0,T], \\
u(0)=0 \end{array}\right.
\end{equation}
has the $L^p-L^2(\R^N)$ maximal regularity.

\section{Maximal regularity for a class of Schr\"{o}dinger operators}

In this section we examine $L^p-L^q$ a priori estimates for the Schr\"odinger type equation (\ref{proSch}). Let again  
$$A(t) = -\Delta + m(t,\cdot)$$
be a  Schr\"odinger operator on $L^2(\R^N, dx)$ ($A(t)$ is defined at the end of  the previous section). We assume that the potential $m$ satisfies
(\ref{a4}) and (\ref{a5}) with some constant $\beta > 1/2.$ We have seen that by  Theorem \ref{monn}  the family $\{-\Delta +m(t,\cdot),\ t\in[0,T]\}$ belongs to $MR(p,2)$. We recall that it means that for each $f\in L^p(0,T;L^2(\R^N))$ there exists a unique $u\in W^{1,p}(0,T;L^2(\R^N))$  with $t\to A(t)u(t)\in L^p(0,T;L^2(\R^N))$ satisfying (\ref{proSch}) in the $L^p(0,T;L^2(\R^N))$ sense. The closed graph theorem implies then that there exists a constant $C>0$ such that 
$$\|u\|_{L^p(0,T;L^2(\R^N))}+\|\partial_t u\|_{L^p(0,T;L^2(\R^N))}+\|Au\|_{L^p(0,T;L^2(\R^N))}\leq C\|f\|_{L^p(0,T;L^2(\R^N))}.$$
By \cite[Theorem 2.1, Proposition 3.11]{CMS}, for every $\lambda>0$, $1\leq p\leq \infty$, the operator
$\lambda+\partial_t+A(t)$ is invertible in $L^p([0,T]\times\R^N)$. Therefore, by adding  a positive constant and  observing that this does not affect the study of the maximal regularity, we can suppose that $\partial_t+A(t)$ is invertible. The maximal $L^p-L^2$ regularity property can be reformulated by saying that the operator $\partial_t(\partial_t-\Delta+m)^{-1}$ is bounded in $L^p(0,T; L^2(\R^N))$. If we prove that such a operator and its adjoint are of weak type $(1,1)$, then the Marcinkiewicz interpolation Theorem  implies that $\{A(t),\ t\in[0,T]\}$ belongs to the class $MR(p,q)$. Here $A(t)$ is a certain realisation of the operator on $L^q(\R^N).$ \\

 From now on, we assume without loss of generality that $\partial_t-\Delta+m$ is invertible on $L^2([0,T] \times \R^N).$  We also assume that
  (\ref{a4}) and (\ref{a5}) with some constant $\beta > 1/2.$ As explained, the operator $\partial_t(\partial_t-\Delta+m)^{-1}$ is bounded on $L^2([0,T] \times \R^N).$ 
 We prove  that $\partial_t(\partial_t-\Delta+m)^{-1}$ is of weak type $(1,1)$.

\begin{prop}  \label{weak}
The operator  $\partial_t(\partial_t-\Delta+m)^{-1}$ is of weak type $(1,1)$.
\end{prop}

{\sc Proof.} Let $f\in L^1([0,T]\times\R^N)\cap L^2(0,T\times \R^N)$ and set  $u :=(\partial_t-\Delta+ m)^{-1}f$ with $u(0,\cdot)=0.$  Note in passing  that 
$u \in  L^1([0,T]\times\R^N)$ (this follows for example from (\ref{kato-cons})  below). \\
Let $h_n:\R\rightarrow \R$ be a sequence of smooth functions such
that $|h_n|\leq C,$ $h_n'(s)\geq 0$ and  $h_n(s)\rightarrow
\sign(s)$ pointwise as  $n\to\infty$. Let $H_n$ be such that $H'_n=h_n$ and
$H_n(0)=0.$ By the Lebesgue convergence
Theorem, we have
\begin{eqnarray*}
\int_{[0,T]\times\R^{N}}\sign(u)\partial_tu&=&\lim_n\int_{[0,T]\times\R^{N}}h_n(u)\partial_tu
=\lim_n\int_{[0,T]\times\R^{N}}\partial_t (H_n(u))\\&=&\lim_n\int_{\R^N}H_n(u(T,\cdot))=\int_{\R^N}\lim_nH_n(u(T,\cdot))\geq 0,
\end{eqnarray*}
and, by dissipativity,
\begin{equation}  \label{dissip}
-\int_{[0,T]\times\R^{N}}\sign(u)\Delta
u
\geq 0.
\end{equation}
Therefore, after multiplying  by $\sign(u)$ both member of 
\begin{equation} \label{eqSchr}
\partial_t u-\Delta u+mu=f,
 \end{equation}  we obtain
$$
\int_{[0,T]\times\R^{N}}m|u|\leq
\int_{[0,T]\times\R^{N}}\sign(u)(\partial_t-\Delta+m)u=\int_{[0,T]\times\R^{N}}f\sign(u)\leq
\int_{[0,T]\times\R^{N}}|f|.
$$
From this and the equation (\ref{eqSchr}) we deduce  that 
\begin{equation} \label{L1bound}
\int_{[0,T]\times\R^{N}}|(\partial_t-\Delta)u|\leq
 2\int_{[0,T]\times\R^{N}}|f|
\end{equation}
(the last inequality is known, it has been proved  by Kato \cite{kato} in the elliptic case. See  also \cite{CMS} for its parabolic version).
Since the operator $\partial_t(\partial_t-\Delta)^{-1}$ is of weak type $(1,1)$ (see for example \cite[Section 5]{hieber-pruss}), it follows that there exists a positive constant $C$ such that for every $\alpha>0$
\begin{eqnarray*}
\mu\{|\partial_t(\partial_t-\Delta+m)^{-1}f|>\alpha\}&=&\mu\{|\partial_t(\partial_t-\Delta)^{-1}(\partial_t-\Delta)(\partial_t-\Delta+m)^{-1}f|>\alpha\}\\&\leq& \frac{C}{\alpha}\|(\partial_t-\Delta)(\partial_t-\Delta+m)^{-1}f\|_{L^1([0,T]\times\R^N)}\leq \frac{C}{\alpha}\|f\|_{L^1([0,T]\times\R^N)}
\end{eqnarray*}
and so $\partial_t(\partial_t-\Delta+m)^{-1}$ is of weak type $(1,1)$.
\qed

By interpolation, from the last proposition we get the maximal $L^p-L^q$ regularity for $1<q\leq 2$. In order to extend the range of admissible values for $q$, we apply a duality argument.

\begin{prop}  \label{weak*}
The adjoint operator  $[\partial_t(\partial_t-\Delta+m)^{-1}]^*$ acting on $C_c^\infty([0,T]\times\R^N)$ is of weak type $(1,1)$.
\end{prop}

One of the  tools  in the proof is a distributional inequality proved by Kato for the
Laplacian (see \cite[Theorem X.2]{simon}). For completeness we provide here a
short proof in the parabolic case and for complex values functions $u$.

\begin{lem}[Parabolic Kato's inequality] \label{kato}
Let $u \in L^1_{loc}([0,T]\times \R^N)$ be such that $(\partial_t-\Delta)u \in
L^1_{loc}([0,T]\times \R^N)$. Define
\begin{eqnarray*} \sign(u)=
\begin{cases}
0&if\ \ \  u(x)=0\\
\overline{u(x)}/|u(x)|&if\ \ \  u(x)\neq 0.
\end{cases}
\end{eqnarray*}
Then $|u|$ satisfies the following distributional inequality
$$
(\partial_t-\Delta)|u|\leq \Rp{[\sign(u)(\partial_t-\Delta)u]}.
$$
\end{lem}
{\sc Proof.}
We first suppose that $u\in C^\infty([0,T]\times \R^N)$. Define
\begin{equation}\label{eq: kato's 1}
    u_\eps(x) :=\sqrt{|u|^2+\eps^2}
\end{equation}
so that $u_\eps\in C^\infty([0,T]\times \R^N).$ Since
\begin{equation}\label{eq: kato's 2}
u_\eps\nabla u_\eps=\Rp[\overline{u}\nabla u].
\end{equation}
and $u_\eps\geq |u|,$ then (\ref{eq: kato's 2}) implies that
\begin{equation}\label{eq: kato's 3}
|\nabla u_\eps|\leq |\overline{u}||u_\eps|^{-1}|\nabla u|\leq |\nabla u|.
\end{equation}
Taking the divergence of (\ref{eq: kato's 2}) we obtain
\begin{equation*}
    u_\eps\Delta u_\eps+|\nabla u_\eps|^2=\Rp(\overline{u}\Delta
    u)+|\nabla u|^2
\end{equation*}
so by (\ref{eq: kato's 3})
\begin{equation}\label{eq: kato's 4}
  \Delta u_\eps\geq \Rp[{\rm sign}_\eps(u)\Delta u],
\end{equation}
where ${\rm sign}_\eps(u) :=\overline{u}/u_\eps.$
Differentiating (\ref{eq: kato's 1}) with respect to $t$ we obtain
\begin{equation}\label{eq: kato's 5}
\partial_t u_\eps=\Rp[{\rm sign}_\eps(u)\partial_t u]
\end{equation}
and combining (\ref{eq: kato's 4}) and (\ref{eq: kato's 5}) yields
\begin{equation}\label{eq: kato's 6}
(\partial_t-\Delta)u_\eps\leq \Rp[{\rm
sign}_\eps(u)(\partial_t-\Delta) u].
\end{equation}
Let now $u\in L^1_{loc}([0,T]\times \R^N)$ be such that $(\Delta-\partial_t)u\in
L^1_{loc}([0,T]\times \R^N)$ and let $\phi_n$ be an approximate identity. Since
$u^n :=u*\phi_n\in C^\infty([0,T]\times \R^N),$ then by (\ref{eq: kato's 6})
\begin{equation}\label{eq: kato's 7}
(\partial_t-\Delta)(u^n)_\eps\leq \Rp[{\rm
sign}_\eps(u^n)(\partial_t-\Delta) u^n].
\end{equation}
Fix $\eps>0$ and let $n\rightarrow \infty.$ Then $u^n\rightarrow
u$ in $L^1_{loc}([0,T]\times \R^N)$ and a.e. (passing to a subsequence, if necessary).
Thus ${\rm sign}_\eps(u^n)\rightarrow {\rm sign}_\eps(u)$
a.e. Since $(\partial_t-\Delta)u^n=\left ((\partial_t-\Delta)u\right
)*\phi_n$ and $(\partial_t-\Delta)u \in L^1_{loc}([0,T]\times \R^N)$, then
$(\partial_t-\Delta)u^n\rightarrow (\partial_t-\Delta)u$ in
$L^1_{loc}([0,T]\times \R^N)$, too. It is now easy to see that ${\rm
sign}_\eps(u^n)(\partial_t-\Delta)u^n$ converges in the sense of
distributions to ${\rm sign}_\eps(u) (\partial_t-\Delta)u.$
Thus, letting $n\rightarrow \infty$ in (\ref{eq: kato's 6}) we
conclude that $$ (\partial_t-\Delta)u_\eps\leq \Rp[{\rm
sign}_\eps(u)(\partial_t-\Delta) u].
$$
Now taking $\eps\rightarrow 0$ we obtain the desired inequality
for $u$, since ${\rm sign}_\eps(u)\rightarrow {\rm sign}(u)$ and
$|{\rm sign}_\eps(u)| \leq 1$.
\qed

We will need also the following simple equality.

\begin{lem}     \label{der-sign}
Let $f,\ g\in C^1([0,T]\times\R^N)$. Then
$$\partial_t(\sign (g)f)=\sign(g)\partial_t f$$
almost everywere in $[0,T]\times\R^N$. Here $\partial_t$ is the distributional derivative with respect to $t$. 
\end{lem}
{\sc Proof.}  Given $\eps>0$, define
    $$g_\eps(x)=\frac{g}{\sqrt{|g|^2+\eps}}$$
so that $g_\eps\in C^1([0,T]\times\R^N)$, $g_\eps\to \sign(g)$ pointwise as $\eps\to 0$ and
$$\partial_t (g_\eps f)=g_\eps \partial_t f+f\left(\frac{\partial_tg(|g|^2+\eps)-|g|^2\partial_t g}{(|g|^2+\eps)^\frac{3}{2}}\right).$$
By letting $\eps$ to $0$, the right hand side in the previous equality converges pointwise to $\sign(g)\partial_t f$ since the second addendum converges obviously to $0$ in the set $\{g\neq 0\}$ and $\partial_t g= 0$ almost everywhere  where $g=0$ (see for example \cite[Lemma 7.7]{GilTru}). The left hand side approachs  $\partial_t(\sign (g)f)$ in the sense of distributions. By observing that, by dominated convergence, the convergence in the right hand side is also in the distributional sense, we deduce the distributional equality
$$\partial_t(\sign (g)f)=\sign(g)\partial_t f$$
and so the claim.\qed

{\sc Proof of Proposition \ref{weak*}.}  Let $g\in C_c^\infty([0,T]\times\R^N)$ and set $u=(\partial_t-\Delta+m)^{-1}g$. Then, by Kato's inequality (Lemma \ref{kato}), 
$$(\partial_t-\Delta)|u|\leq (\partial_t-\Delta+m)|u|\leq \sign (u)(\partial_t-\Delta+m)u=\sign (u) g$$ and, since $(\partial_t-\Delta)^{-1}$ is a positive operator, $|u|\leq (\partial_t-\Delta)^{-1}(\sign (u) g)$. We have from the definition of $u$,
\begin{equation}  \label{kato-cons}
|(\partial_t-\Delta+m)^{-1}g|\leq (\partial_t-\Delta)^{-1}(\sign (u) g)\leq (\partial_t-\Delta)^{-1}|g|.
\end{equation}
Let us consider now the adjoint operator $[\partial_t(\partial_t-\Delta+m)^{-1}]^*=[(\partial_t-\Delta+m)^{-1}]^*\partial_t^*$.
By (\ref{kato-cons}) we have that 
$|[(\partial_t-\Delta+m)^{-1}]^*g|\leq [(\partial_t-\Delta)^{-1}]^*|g|$ for all $g\in C_c^\infty([0,T]\times\R^N)$. Therefore, if $f\in C_c^\infty([0,T]\times\R^N)$,  by choosing $g=\partial_t^* f$, after an approximation procedure as in Lemma \ref{der-sign} and by recalling that $\partial_t^*=-\partial_t$ where it is defined and that the operator $[(\partial_t-\Delta)^{-1}]^*$ is continuous 
on $L^2([0,T]\times\R^N)$, we deduce that 
$$[(\partial_t-\Delta)^{-1}]^*|\partial_t^* f|= [\partial_t(\partial_t-\Delta)^{-1}]^*(sign(\partial_t^* f) f).$$
Therefore,
\begin{eqnarray*}
|[\partial_t(\partial_t-\Delta+m)^{-1}]^*f|&=&|[(\partial_t-\Delta+m)^{-1}]^*\partial_t^*f|\\
&\leq&  [(\partial_t-\Delta)^{-1}]^*|\partial_t^* f|\\
&=&[\partial_t(\partial_t-\Delta)^{-1}]^*(sign(\partial_t^* f) f).
\end{eqnarray*}
Then, since $[\partial_t(\partial_t-\Delta)^{-1}]^*$ is of weak type (1,1) (see \cite[Section 5]{hieber-pruss}), there exists a positive constant $C$ such that for every $f\in C_c^\infty([0,T]\times\R^N),\ \alpha>0,$ 
\begin{eqnarray*}
\mu\{|[\partial_t(\partial_t-\Delta+m)^{-1}]^*f|>\alpha\}&\leq&\mu\{|[\partial_t(\partial_t-\Delta)^{-1}]^*(sign(\partial_t^* f) f)|>\alpha\}\\&\leq& \frac{C}{\alpha}\|sign(\partial_t^* f) f\|_{L^1([0,T]\times \R^N)}=\frac{C}{\alpha}\| f\|_{L^1([0,T]\times \R^N)}
\end{eqnarray*} 
and the proof is complete.\qed

\begin{os} \label{appr}
By  Proposition \ref{weak*} and by approximation we deduce that the adjoint operator  $[\partial_t(\partial_t-\Delta+m)^{-1}]^*$ is of weak type $(1,1)$  on $L^1([0,T]\times\R^N)\cap L^2([0,T]\times \R^N)$. Indeed it is sufficient to approximate a given function $f\in L^1([0,T]\times\R^N)\cap L^2([0,T]\times \R^N)$ with smooth functions $(f_n)\subset C_c^\infty([0,T]\times\R^N)\cap L^2(0,T\times \R^N)$  in the $L^1([0,T]\times\R^N)\cap L^2([0,T]\times \R^N)$ norm. Then the convergence of the set's measures 
$\mu\{\{|[\partial_t(\partial_t-\Delta+m)^{-1}]^*f_n|>\alpha\}\}$  follows by the $L^2([0,T]\times\R^N)$ boundedness of the operator $[\partial_t(\partial_t-\Delta+m)^{-1}]^*$.
\end{os}

We can now state the main theorem of this section.  Before  that we clarify what we mean by $-\Delta+m$ on $L^p([0,T]\times \R^N).$  If $p = 2$ this operator is constructed 
by sesquilinear forms (see Example 2 at  the end of the  previous section). As explained previously, 
$\la I + \partial_t -\Delta+m$ is invertible on $L^2([0,T]\times \R^N)$ for all $\la > 0$ and $(\la I + \partial_t -\Delta+m)^{-1}$ defines  a bounded operator on $L^p([0,T]\times \R^N)$ (see
(\ref{kato-cons})). By a simple density argument, it satisfies the resolvent  equation and hence it is the resolvent of a certain closed operator on $L^p([0,T]\times \R^N)$. We denote,  as in the case $p = 2,$ this operator by $\partial_t -\Delta+m.$ Hence, we have a realisation of  the operator $-\Delta+m$ on $L^p([0,T]\times \R^N)$ such that $\la I + \partial_t -\Delta+m$ is invertible and  $(\la I + \partial_t -\Delta+m)^{-1}$ coincides  on $L^p([0,T]\times \R^N) \cap L^2([0,T]\times \R^N)$ with the resolvent of the starting operator on $L^2([0,T]\times \R^N).$ For this realisation, the Cauchy problem (\ref{proSch}) has a unique solution in the $L^p([0, T]\times\R^N)$ sense.

By Propositions \ref{weak}, \ref{weak*},  Theorem \ref{monn} and   the Marcinkiewicz interpolation Theorem, we  deduce that the operator $\partial_t(\partial_t-\Delta+m)^{-1}$ is bounded in $L^p(0,T; L^q(\R^N))$ for all $1<p,\ q<\infty$. Then we have proved the following result.

\begin{teo} \label{thm} 
Let $0\leq m(t,x)\in L^1_{loc}([0,T]\times \R^N)$. 
As before, we assume that there exists a non-negative potential $W \in L^1_{loc}(\R^N, dx)$ such that $m$ satisfies the following properties  (in which $c_1, c_2$    are positive constants and $\beta > 1/2$)
$$
c_1 W(x) \le m(t,x) \le c_2 W(x) \ (a.e.\  x \in \R^N) \ \mbox{and all } \ t \in [0,T], 
$$
$$
 \vert m(t,x) - m(s,x) \vert  \le c_2 W(x) \vert t - s \vert^\beta  \ (a.e.\ x \in \R^N) \ \mbox{and all } \ t, s \in [0,T].$$
Then, for $1<p,\ q<\infty$, the family $\{A(t)=-\Delta+m(t,\cdot),\ t\in [0,T]\}$ belongs to the class $MR(p,q)$.
\end{teo}

\section{Extension to other operators} 
In the previous section, we gave in details the proof of the $L^p-L^q(\R^N)$ maximal  regularity   of the operator $-\Delta + m(t,\cdot)$. 
See Theorem \ref{thm}.  Here we explain how to extend this result to  more general situations. Since the proofs are similar to those in the previous section 
we shall not rewrite all the details but mention the main ingredients. 
 
 Let $(X, \mu, \rho)$ be a metric measured space.  We denote by $v(x, r)$ the volume of the ball of center $x$ and radius $r,$ that is
 $$ v(x, r) := \mu \left( B(x, r)\right) := \mu\left( \{ y \in X, \rho(x,y) < r \} \right).$$
 We assume that $v(x, r) < \infty$ for all $x \in X$ and $r > 0$ and that $X$  satisfies the doubling condition
 \begin{equation}\label{doubl1}
 v(x,2r) \le C_{0} v(x, r) \ \forall x \in X, \ r > 0,
 \end{equation}
 where $C_{0}$ is a positive constant (independent of $x$ and $r$). 
 
 Let now $\Om$ be an open subset of $X$ and consider $A : D(A)\subseteq L^2(\Om, \mu) \to L^2(\Om, \mu)$ be a densely defined linear operator.  We assume for simplicity that 
 $A$ is a non-negative self-adjoint operator and denote by $(e^{-tA})_{t\ge 0}$ its associated semigroup on  $ L^2(\Om, \mu).$ 
 We assume that   $(e^{-tA})_{t\ge 0}$ is a  sub-markovian semigroup. This implies in particular that  $(e^{-tA})_{t\ge 0}$ acts as a strongly continuous semigroup on $L^p(\Om, \mu)$ for $1 \le p < \infty.$ For simplicity, we keep the same notation as in $L^2(\Om, \mu)$ and  denote by 
 $-A$ the corresponding generator on $L^p(\Om, \mu)$. 
  Finally, we assume that $e^{-tA}$ is given by a kernel $p(t, x, y)$ 
 (called the heat kernel of $A$) which 
 satisfies the global  Gaussian upper bound
 \begin{equation}\label{UE} 
 \vert p(t, x, y) \vert \le  \frac{C}{ \sqrt{v(x, t^{1/m}) v(y, 
 t^{1/m})}} \exp \left\{-c \frac{\rho(x,y)^{m/(m-1)}}{t^{1/(m-1)}} \right\} 
 \end{equation}
for all $t > 0$ and $\mu$-a.e. $x, y \in \Om.$  Here $C,  c$ and $m$ are positive constants and $m \ge 2$. \\
The above assumptions are satisfied for a wide class of operators including divergence form uniformly elliptic operators on domains (with Dirichlet boundary conditions), Schr\"odinger operators and also Laplace-Beltrami operators on some Riemannian manifolds, see \cite{davies} or \cite{ouhabaz}. 

Given a non-negative potential   $m : [0,T] \times \Om \to [0, \infty]$, we construct  self-adjoint operators
$$A(t) := A + m(t, \cdot)$$
in the same way as $-\Delta + m(t, \cdot)$ in Section \ref{sec3}. Each $A(t)$ is the associated operator with the form
$$\fra(t,u,v) := \fra(u,v) + \int_\Om m(t,\cdot) u v dx$$
 where $\fra$ is the form of the self-adjoint operator $A$. We assume that $\fra(t, \cdot, \cdot)$ (defined on the intersection domain) is densely defined and closed and 
 assume that (\ref{a4}) and (\ref{a5}) are satisfied with some constant $\beta > 1/2$ (for $\mu-a.e. x \in \Om$). We conclude  from Theorem \ref{monn} that $L^p-L^2(\Om, \mu)$
 maximal regularity holds for the family $A(t), t \in [0,T].$ In order to obtain $L^p-L^q(\Om, \mu)$ estimates we proceed as in the previous section. 
 
 Consider first the case $1 < q < 2. $ 
  The sub-Markovian assumption of the semigroup implies in particular that the corresponding generator on $L^1(\Om, \mu)$ is accretive. That is 
 $$\int_\Om A u\cdot sign u\, d\mu \ge 0$$
 for all $u$ in  the $L^1-$domain. This implies (\ref{L1bound}) with $A$ in place of $-\Delta.$ 
The doubling condition (\ref{doubl1}) and the Gaussian upper bound (\ref{UE}) imply that the operator 
$\partial_t (\partial_t + A)^{-1}$ is of weak type (1,1). See \cite{hieber-pruss} and \cite{CoDu}. Using this, the proof of Proposition \ref{weak} works without any modification. 
 
 In order to treat the case $2 < q < \infty,$  we proceed as in Proposition \ref{weak*}.  For  any positive constant $\lambda,$ the  operator $(\lambda I + \partial_t + A +m(t,\cdot) )^{-1}$ defines a  bounded operator on  $L^q([0,T] \times \Om)$
 for $1 \le q \le \infty.$  Indeed, it is bounded on $L^2([0,T] \times \Om)$ (by the $L^2-L^2$ maximal regularity) and since the semigroup $(e^{-tA})_{t\ge 0}$ is positive (since it is sub-Markovian) and $m$ is non-negative then the pointwise inequality 
 \begin{equation}\label{domi}
 0 \le  (\lambda I + \partial_t + A +m(t,\cdot) )^{-1} g  \le  (\lambda I + \partial_t + A  )^{-1}g 
 \end{equation}
holds for a.e. $(t,x) \in [0,T] \times \Om$ and all non-negative function $g$ (this can be checked by applying the positivity and domination criteria in \cite[Chapter 2]{ouhabaz} to the operators
$\partial_t + A + m(t,\cdot)$ and $\partial_t + A$). This means that (\ref{kato-cons}) holds in this context. The rest of the proof of Proposition \ref{weak*} does not change. Therefore we obtain 
Theorem \ref{thm} for $A(t) = A + m(t, \cdot).$

\end{document}